\documentclass[12pt]{article}
\usepackage{amscd,amsthm}
\usepackage{latexsym}
\usepackage{amsmath}
\usepackage{amssymb}
\numberwithin{equation}{section}

\def\bbR{{\Bbb R}}
\def\m{{\cal M}}
\def\k{{\cal K}}
\def\p{{\cal P}}
\def\u{{\cal U}}

\def\const{{\rm const}}
\def\Trace{{\rm Trace}}

\begin{document}
\title{Central Limit Theorem for Local Linear Statistics in Classical
Compact Groups and Related Combinatorial Identities}
\author{Alexander Soshnikov\\California Institute of Technology\\Department of
Mathematics\\Sloan 253-37\\Pasadena, CA  91125 USA\\and\\ University of
California, Davis\\Department of Mathematics\\Davis, CA  95616, USA}
\date{}
\maketitle
\begin{abstract}
We discuss CLT for the global and local linear statistics of random
matrices from classical compact groups.  The main part of our proofs
are certain combinatorial identities much in the spirit of  works
by Kac and Spohn.
\end{abstract}

\section{Introduction}

Let $M$ be a unitary matrix chosen at random with respect to the Haar 
measure on the unitary group $U(n)$.  We denote the eigenvalues of $M$
by $\{\exp (i\cdot\theta_j)\}^n_{j=1}$, where $-\pi\leq \theta_1,\theta_2,
\ldots ,\theta_n<\pi$.  The joint distribution of the eigenvalues
(called the Weyl measure) is absolutely continuous with respect to the
Lebegue measure $\prod^n_{j=1}d\theta_j$ on the $n$-dimensional tori and its
density is given by
\begin{equation}
P_{U(n)}(\theta_1,\ldots ,\theta_n)=\frac{1}{(2\pi )^n\cdot n!}
\cdot\prod_{1\leq j<k\leq n}\vert\exp (i\cdot\theta_j)-\exp (i\cdot\theta_k)
\vert^2
\end{equation}
(see [We]).  Throughout the paper we will be interested in the global and
local linear statistics
\begin{equation}
S_n(f)=\sum^n_{j=1}f(\theta_j),
\end{equation}

\begin{equation}
S_n(g(L_n\cdot ))=\sum^n_{j=1}g(L_n\cdot\theta_j),
\end{equation}

$$L_n\rightarrow\infty ,\ \frac{L_n}{n}\rightarrow 0.$$
The optimal conditions on $f,g$ for our purposes are
\begin{equation}
\sum^\infty_{k=-\infty}\vert\hat f(k)\vert^2\cdot\vert k\vert <\infty ,
\end{equation}

\begin{equation}
\int^\infty_{-\infty}\vert\hat g(t)\vert^2\cdot\vert t\vert dt <\infty ,
\end{equation}
where
\begin{align*}
f(x)&=\sum^\infty_{k=-\infty}\hat f(k)\cdot e^{ikx},\\
g(x)&=\frac{1}{\sqrt{2\pi}}\cdot\int^\infty_{-\infty}\hat g(t)\cdot
e^{itx}dt
\end{align*}
However in order to simplify the exposition we will always assume that
$f$ has a continuous derivative on a unit circle ( $f\in C^1(S^1)$ )
and $g$ is a Schwartz function ( $g\in f(J(\bbR^1))$ ).

Let us denote by $E_n$ the mathematical expectation with respect to Haar
measure.  We start with the formulation of the result which is essentially 
due to C. Andr\'eief ([A], for a modern day reference see [TW]and also [Dy].
\medskip

\noindent{\bf Proposition}
\begin{equation}
E_n \exp (t S_n(f))-\det (Id+(e^{t f}-1) K_n)=\det (Id+(e^{t
f}-1) Q_n),
\end{equation}
where $(e^{t f}-1)$ is a multiplicaiton operator and $K_n, Q_n:L^2(S^1) \to L^2(S^1)$
are the integral operators with the kernels
\begin{equation}
K_n(x,y)=\frac{1}{2\pi}\frac{\sin\left (\frac{n}{2} (x-y)\right )}
{\sin\left (\frac{x-y}{2}\right )},
\end{equation}

\begin{equation}
Q_n(x,y)=\sum^{n-1}_{j=0}\frac{1}{\sqrt{2\pi}}e^{ijx}\frac{1}{\sqrt
{2\pi}}e^{-i j y}
\end{equation}

\noindent{\bf Remark 1}.  $K_n, Q_n$ are unitary equivalent to each other
and are the operators of a finite rank.  In particular, $Q_n$ is just a projection
operator on the first $n$ harmonic functions of the unit circle.

One of the ingredients of the proof of the proposition 
is the following chain of the equalities
\begin{align}
\begin{split}
p_{U(n)}(\theta_1,\ldots ,\theta_n)&=\frac{1}{n!}\cdot\det (e^{i\cdot 
(j-1)\cdot \theta_k})_{1\leq j,k\leq n}\cdot\det (e^{-i\cdot (j-1)\cdot
\theta_k})_{1\leq j,k\leq n}\\
&=\frac{1}{n!}\det\biggl (Q_n(\theta_j ,\theta_k)\biggr )_{1\leq j,k\leq n}\\
&=\frac{1}{n!}\det\biggl (K_n(\theta_j ,\theta_k)\biggr )_{1\leq j,k\leq n}
\end{split}
\end{align}
Remark 1 allows us to rewrite the Fredholm determinants in (1.6) as the
Toeplitz determinant with the symbol $\exp \bigl (t\cdot f(\cdot )\bigr )$:
\begin{align}
\begin{split}
E_n\exp\biggl (t\sum^n_{j=1}f(\theta_j)\biggr )&= D_{n-1}\biggl (
\exp (t\cdot f)\biggr )\\
&=\det \biggl (\frac{1}{2\pi}\int^{2\pi}_0\exp\bigl (t f(x)\bigr )
\cdot\exp\bigl (i (j-k)x\bigr )dx\biggr )_{1\leq j,k\leq n}
\end{split}
\end{align}
The asymptotics of (1.10) for large $n$ is given by the Strong Szego Limit
Theorem:
\begin{align}
\begin{split}
&D_{n-1}\biggl (\exp (t\cdot f)\biggr )=\\
&\qquad \exp\biggl (t n \hat f(0)+\frac{1}{2}t^2
\sum^{+\infty}_{-\infty}\vert k\vert \vert\hat f(k)\vert^2+\bar 0
(1)\biggr )
\end{split}
\end{align}
(see [Sz] and  [K], [H], [De], [F-H], [G-I], [Wid1], [Wid2], [McC-W], [Ba-W], [Jo1], [Bo], [Bo-S], [Me],
[So2], [Wie], [D] for further developments.)

In probabilistic terms (1.11) claims that $ES_n(f)=\tfrac{n}{2\pi}\cdot
\int^\pi_{-\pi}f(\theta )d\theta +\bar 0(1)$ (actually the remainder term
is  zero), and the centralized random variable $\sum^n_{j=1}
f(\theta_j)=E_n\sum^n_{j=1}f(\theta_ j)$ converges in distribution to the
normal law $N(0,\sum^\infty_{-\infty}\vert k\vert \vert\hat f(k)\vert^2$).

Our first goal is to establish a similar result for the local linear
statistics.
\medskip

\noindent{\bf Theorem 1}.  {\it Let $g\in J(\bbR^{1}),\ L_n\rightarrow
+\infty,\ \tfrac{L_n}{n}\rightarrow 0$.  Then 
$E_n\sum^n_{j=1}g(L_n\cdot\theta_j )
=\tfrac{n}{2\pi\cdot L_n}\cdot\int^\infty_{-\infty} g(x)dx$, and the
centralized random variable $\sum^n_{j=1}(g(L_n\cdot\theta_j)-E
\sum^n_{j=1}g(L_n \theta_j)$ converges in distribution to the normal
law $N(0,\tfrac{1}{2\pi}\cdot\int^{+\infty}_{-\infty}\vert\hat g(t)\vert^2
\vert t\vert dt)$.}
\medskip

We give a combinatorial proof which holds  both in the local and global cases.  In some sense
our approach is close to the  heuristic arguments in
[I-D].  We start with
\medskip

\noindent{\bf Lemma 1}.  {\it Let $C_{\ell ,n}(f)$ be the $\ell$-th cumulant
of $S_n(f)$.  Then
\begin{align}
\begin{split}
&\vert C_{\ell ,n}(f)-\sum_{k_1+\ldots +k_\ell =0}\hat f(k_1)\cdot\ldots\cdot
\hat f(k_\ell )\cdot\sum^\ell_{m=1}\frac{(-1)^{m-1}}{m}\cdot\\
&\qquad\sum_{\stackrel{\text{\footnotesize
$\ell_1+\ldots +\ell_m=\ell ,$}}{\ell_1\geq 1,
\ldots ,\ell_m
\geq 1}}\frac{\ell !}{\ell_1!\cdot\ldots\cdot\ell_m!}\cdot \biggl (n-{\rm max}
(0,\sum^{\ell_1}_{i=1}k_i , \sum^{\ell_1+\ell_2}_{i=1}k_i,\ldots ,\\
&\qquad \sum^{\ell_1+\ldots +\ell_{m-1}}_{i=1}k_i)-{\rm max}(0,
\sum^{\ell_1}_{i=1}(-k_i),\sum^{\ell_1+\ell_2}_{i=1}(-k_i),\ldots ,\\
&\qquad \sum^{\ell_1+\ldots +\ell_{m-1}}_{i-1}(-k_i)) \biggr )\vert\leq
\const_\ell
\cdot \sum_{\stackrel{\text{\footnotesize $k_1+\ldots +k_\ell =0$}}
{\vert k_1\vert +\ldots
+\vert k_\ell\vert >n}}\vert k_1 \vert \vert \hat f(k_1)\vert\cdot\ldots\cdot\vert\hat f
(k_\ell )\vert
\end{split}
\end{align}
}
\medskip
\break
\noindent{\bf Remark 2}  One can see that for sufficiently smooth $f$ 
the r.h.s. of (1.12) goes to zero
as $n\rightarrow\infty $.
\medskip

\noindent{\bf Remark 3}  An analogous result to lemma 1 was established
in [Spo] for the determinantal random point field with the sine kernel
(see also Remark 4 below).

The proof of Lemma 1 will be given in \S 2.  At this state we observe
that it implies
\medskip

\noindent{\bf Lemma 2}  {\it The limit of $C_{\ell ,n}(f), \ell >1$
exists as $n\rightarrow \infty$ and is equal to $\sum_{k_1+\ldots +k_\ell =0}
\hat f(k_1)\cdot\ldots\cdot\hat f(k_\ell )\cdot (G(k_1,\ldots ,k_\ell)
+G(-k_1,\ldots ,-k_\ell ))$, where $G$ is the piece-wise linear continuous
function defined by
\begin{align}
\begin{split}
&G(k_1,\ldots ,k_\ell ):=\sum_{\sigma\in S_\ell}\sum^\ell_{m=1}\frac
{(-1)^m}{m}\cdot
\sum_{\stackrel{\text{\footnotesize$\ell_1+\ldots +\ell_m=\ell ,$}}
{\ell_1\geq 1,\ldots ,
\ell_m\geq 1}}\frac{1}{\ell_1!\cdot\ldots\cdot\ell_m!}\cdot\\
&\qquad {\rm max} \left (0, \sum^{\ell_1}_{i=1}k_{\sigma (i)},
\sum^{\ell_1+\ell_2}_{i=1}k_{\sigma (i)},\ldots ,\sum^{\ell_1+\ldots
+\ell_{m-1}}_{i=1}k_{\sigma (i)}\right ).
\end{split}
\end{align}
}
\medskip

\noindent{\bf Proof of Lemma 2}  After opening the brackets in (1.12) we observe
that the coefficient in front of $n$ is equal to
\begin{equation}
\sum^\ell_{m=1}\sum_{\stackrel{\text{\footnotesize $\ell_1+\ldots +\ell_m=
\ell ,$}}{\ell_1\geq 1,
i=1,\ldots ,m}}\frac{(-1)^{m-1}}{m}\frac{\ell !}{\ell_1!\ldots\ell_m!}=
\begin{cases}
1, &\ell =1\\0,&\ell >1
\end{cases}
\end{equation}
Indeed, the generating function of these coefficients is equal to
$$\log\biggl (1+(e^z-1)\biggr )=z.$$ 
\qed

Now CLT for $\sum^n_{j=1} f(\theta_j)$ follows from

\begin{center}
{\bf Main Combinatorial Lemma}
\end{center}

{\it Let $k_1,\ldots ,k_\ell$ be arbitrary real numbers such that their sum equals
zero.  Let $G(k_1,\ldots ,k_\ell )$ be defined as in (1.13).  Then
$$G(k_1,\ldots ,k_\ell )=\begin{cases}\vert k_1\vert =\vert k_2\vert &
\text{if }\ell=2\\0 &\text{if }\ell >2\end{cases}.
$$
}
We will prove the lemma in \S 3.
\medskip

\noindent{\bf Remark 4}  A similar combinatorial lemma was stated by
Spohn in [Spo].  He studied a time-dependent motion of a system of infinite
number of particles governed by the equations
$$d\lambda_j(t)=\sum_{i\neq j}\frac{1}{\lambda_i-\lambda_j}dt+db_j(t),$$
where $\{b_j(t)\}^{+\infty}_{j=-\infty}$-independent standard brownian 
motions, and the initial distribution of particles is given by determinantal
random point field with the sine kernel $\tfrac{\sin\pi (x-y)}{\pi (x-y)}$.
However, no correct proof of the combinatorial result was given there.  For
completeness we give a proof of Spohn's lemma independently from the proof
of our Main Combinatorial Lemma  in \S 3.

Assuming the combinatorial part is done we can quickly finish the proof
of Theorem 1.  The formula for the mathematical expectation is trivial.
Rewriting (1.12) for the higher cumulants of $\sum^n_{j=1}g(L_n\cdot
\theta_j)$ we see that the limit of the $\ell$-th cumulant is given by 
$$(2\pi )^{-\frac{\ell}{2}}\cdot\int\hat g(t_1)\cdot\ldots\cdot \hat
g( t_\ell )\cdot\biggl (G(t_1,\ldots ,t_\ell)+G(-t_1,\ldots ,-t_\ell )
\biggr )dt_1\ldots dt_\ell$$
where the integral is over the hyperplane $t_1+\ldots +t_\ell =0$.

Theorem 1 is proven.\hfill$\Box$
\medskip

\noindent{\bf Remark 5}  Our method also gives an elementary combinatorial
proof of Szeg\"o theorem ((1.11)) for $f\in C^1(S^1)$ and sufficiently
small complex $t$. It is different from the one suggested by Kac in [K] 
where  the Taylor expansion of
$D_n (1-tg)$ \  as a function of t was calculated and 
then a so- called Kac-Spitzer combinatorial lemma  was employed to confirm
(1.11).
\medskip

\noindent{\bf Remark 6}  Results similar to Theorem 1 have been established
for other random matrix models in [Spo], [Jo3], [KKP], [Ba], [B-F],
[SSo1], [SSo2], [BM-K].

The rest of the paper is organized as follows.  We prove Lemma 1 in \S 2
and Main Combinatorial Lemma in \S 3.  The result analogous to Theorem 1
for orthogonal and symplectic groups is established in \S 4.

The author would like to thank Ya. Sinai, P.Diaconis, K. Johanson and A. Khorunzhy for
useful discussions. The work was  partially supported by the Euler stipend from the German
Mathematical Society.

\section{Proof of Lemma 1}

We start with calculating the moments of $S_n(f)$.  Le us remember that 
$k$-point correlation function of the eigenvalues of random unitary 
matrix is given by

\begin{align}
\begin{split}
\rho_{n,k}(\theta_1,\ldots ,\theta_k)&=\frac{n!}{(n-k)!}\int_{T^{n-k}}
p_{U(n)}(\theta_1,\ldots ,\theta_n)d\theta_{k+1}\ldots d\theta_n\\
&=\det\biggl (K_n(\theta_i,\theta_j)\biggr )_{1\leq i,j\leq k}=\det
\biggl (Q_n(\theta_i,\theta_j)\biggr )_{1\leq i,j\leq k}
\end{split}
\end{align}
The $N$-th moment of $S_n(f)$ is equal to
$$E_n\left (\sum^n_{i_1=1}f(\theta_{i_1})\cdot\ldots\cdot\sum^n_{i_N=1}
f(\theta_{i_N})\right ),$$
where the indices $i_1,\ldots ,i_N$ range independently from 1 to $n$,
and  in particular  can coincide.  Let $\m =\{M_1,\ldots , M_r\}$ be
a partition of the set $\{1,2,\ldots ,N\}$ into subsets determined by
coinciding indices among $i_1,\ldots ,i_N: M_1=\{j^{(1)}_1,\ldots ,
j^{(1)}_{s_1}\},\ldots , M_r=\{j^{(r)}_1,\ldots ,j^{(r)}_{s_r}\},
\sqcup^r_{i=1}M_i=\{1,2,\ldots ,N\},\break s_i=\vert M_i\vert , 
i=1,\ldots r$.  Then
\begin{equation}
E_n\biggl (S_n(f)\biggr )^N=\sum_{\stackrel{\text{\footnotesize over all}}
{\rm partitions\ \m}} E_n\sum_{\ell_1\neq \ell_2\neq\ldots\neq\ell_r}
f^{s_1}(\theta_{\ell_1})\cdot\ldots\cdot f^{s_r}(\theta_{\ell_r})
\end{equation}
Let us consider a typical term in (2.2) corresponding to a partition $\m$.
\begin{align}
\begin{split}
&E_n\sum_{\ell_1\neq\ldots\neq\ell_r}f^{s_1}(\theta_{\ell_1})\cdot\ldots\cdot
f^{s_r}(\theta_{\ell_r})=\\
&\qquad \int_{T^r}f^{s_1}(x_1)\cdot\ldots\cdot f^{s_r}(x_r)
\cdot\rho_{n,r}(x_1,\ldots ,x_r)dx_1\ldots dx_r
\end{split}
\end{align}
By definition of the determinant and (2.1)
$$\rho_{n,r}(x_1,\ldots ,x_r)=\sum_{\sigma\in S_r}(-1)^\sigma\prod^r_{i=1}
Q_n(x_i,x_{\sigma (i)}).$$
Writing the permutation $\sigma\in S_r$ as a product of cyclic permutations
we have
\begin{align}
\begin{split}
&\rho_{n,r}(x_1,\ldots ,x_r)=\sum_{\stackrel{\text{\footnotesize
over partitions}}{\k\text
{ of }\{1,\ldots r\}}}\biggl (\prod^q_{\alpha =1}\bigl ((-1)^{p_\alpha -1}
\cdot\\
&\qquad \sum_{\stackrel{\text{\footnotesize over all cyclic}}
{\text{permutations of }K_\alpha}}
\prod^{p_\alpha}_{j=1}Q_n(x_{t^{(\alpha )}_j},
x_{\sigma (t^{(\alpha )}_j})\bigr )\biggr )
\end{split}
\end{align}
where $\{1,\ldots ,r\}=\sqcup^q_1K_\alpha, K_\alpha =\{t^{(\alpha )}_1,
\ldots ,t^{(\alpha )}_{p_\alpha}\},\alpha =1,\ldots ,q, p_\alpha =\vert 
K_\alpha \vert$.  Substituting (2.4) into (2.3) we arrive at the expression
that has the following form :
$$
\sum_{\stackrel{\text{\footnotesize over partitions }}
{\m =\{M_1,\ldots , M_r\}\text{ of }\{1,\ldots , N\}}}\qquad
\sum_{\stackrel{\text{\footnotesize over partitions }}
{\k =\{K_1,\ldots , K_q\}\text{ of }\{1,\ldots , r\}}}\cdots .
$$
To interchange the order of summation we construct a new partition $\p =
\{P_1,\ldots , P_q\}$ of $\{1,2,\ldots N\}$ as follows:  $P_i=\sqcup_{j\in
K_i}M_j, i=1,\ldots , q$.  Then $\{M_j\}_{j\in K_i}$ gives a partition of
$P_i$ that we denote by $\p_i$.  We have
\begin{align}
\begin{split}
&E_n\bigl (S_n(f)\bigr )^N=\sum_{\stackrel{\text{\footnotesize over
partitions}}{\p =\{P_1,\ldots ,P_q\}\text{ of }\{1,\ldots ,N\}}}\biggl (
\prod^q_{i=1}\bigl (\sum_{\stackrel{\text{\footnotesize over partitions}}
{\p_i\text{ of }P_i :  \p_i=\{P_{i,1},\ldots ,P_{i,t_i}\}}}\\
&\qquad \int_{T^{t_i}} f^{\vert P_{i,1}\vert}(x_1)\cdot\ldots\cdot
f^{\vert P_{i,t_i}\vert}(x_{t_i})(-1)^{t_i-1}\cdot\sum_{\stackrel
{\text{\footnotesize over cyclic}}{\text{permutations }\sigma\in S_{t_i}}}\\
&\qquad \prod^{t_i}_{j=1} Q_n(x_j,x_{\sigma (j)})dx_1\ldots dx_{t_i}
\bigr )\biggr ).
\end{split}
\end{align}
We remind  that the moments are expressed in terms of cumulants as
$$m_N=\sum_{\stackrel{\text{\footnotesize over partitions}}{ \p=\{P_1,
\ldots ,P_k\}}}C_{\vert P_1\vert}\cdot\ldots\cdot C_{\vert P_k\vert}.$$
Comparing the last formula with (2.5) we arrive at
\begin{align}
\begin{split}
C_{\ell ,n}(f)&=\sum_{\stackrel{\text{\footnotesize partitions}}
{ \p=\{R_1,\ldots , R_m\}\text{ of }\{1,\ldots ,\ell \}}}\int_{T^m}
f^{\vert R_1\vert}(x_1)\cdot\ldots\cdot f^{\vert R_m\vert}(x_m)\cdot\\
&\qquad\qquad (-1)^{m-1}\cdot\sum_{\stackrel{\text{\footnotesize cyclic permutations}}
{\sigma\in S_n}}\prod^m_{j=1}Q_n(x_j,x_{\sigma (j)})\\
&=\sum^\ell_{m=1}\quad \sum_{\stackrel{\text{\footnotesize over ordered
collections}}{(\ell_1,\ldots ,\ell_m): \sum^m_1\ell_i=\ell, \ell_i\geq 1}}
(-1)^{m-1}\frac{\ell !}{\ell_1!\ldots\ell_m!}\frac{1}{m!}\cdot \\ \int_{T^m}
&\qquad\qquad f^{\ell_1}(x_1)\cdot\ldots\cdot f^{\ell_m}(x_m)\cdot m!\cdot
\frac{1}{m}\cdot\\
&\qquad\qquad \prod^m_{j=1} Q_n(x_j,x_{j+1})dx_1\ldots dx_m\\
&=\sum^\ell_{m=1}\quad\sum_{\stackrel{\text{\footnotesize $(\ell_1,
\ldots ,\ell_m): \ell_1+\ldots
+\ell_m=\ell ,$}}{\ell_i\geq 1, i=1, \ldots , m}}
\end{split}
\end{align}
\begin{align}
\begin{split}
&\frac{(-1)^{m-1}}{m}\cdot\frac{\ell !}{\ell_1!\cdot\ldots\cdot\ell_m!}
\cdot \int_{T^m}f^{\ell_1}(x_1)\cdot\ldots\cdot f^{\ell_m}(x_m)\cdot\\
&\qquad \prod^m_{j=1} Q_n(x_j,x_{j+1})dx_1\ldots dx_m.
\end{split}
\end{align}
Since $Q_n(x,y)=\sum^{n-1}_{j=0} e^{-ij(x-y)}$ we can rewrite (2.7) as
\begin{align*}
\begin{split}
&C_{\ell ,n}(f)=\sum^\ell_{m=1}\sum_{\stackrel{\text{\footnotesize
$(\ell_1,\ldots ,\ell_m):$}}
{\ell_1+\ldots +\ell_m=\ell ,\ \ell_i\geq 1}}\frac{(-1)^{m-1}}{m}\cdot
\frac{\ell !}{\ell_1!\ldots\ell_m!}\\
&\qquad\sum^{n-1}_{s_1=0}\ldots
\sum^{n-1}_{s_m=0}\widehat{f^{\ell_1}}(-s_m+s_1)\cdot\widehat{f^{\ell_2}}
(-s_1+s_2)\cdot\ldots\cdot\widehat{f^{\ell_m}} (-s_{m-1}+s_m).
\end{split}
\end{align*}
Writing down the Fourier coefficients of the powers of $f$ as the convolutions
of the Fourier coefficients of $f$
\begin{align*}
\begin{split}
\widehat{f^{\ell_1}}(-s_m+s_1)&=\sum_{\stackrel{\text{\footnotesize 
$(k_1,\ldots ,k_{\ell_1}):$}}
{k_1+\ldots +k_{\ell_1}=s_1-s_m}}\hat f(k_1)\cdot\ldots\cdot \hat 
f(k_{\ell_1}),\\
\widehat{f^{\ell_2}}(-s_1+s_2)&=\sum_{\stackrel{\text{\footnotesize 
$(k_{\ell_1+1},\ldots ,k_{\ell_2}):$}}
{k_{\ell_1+1}+\ldots +k_{\ell_2}=s_1-s_2}}\hat f(k_{\ell_1+1})
\cdot\ldots\cdot \hat f(k_{\ell_2}),\cdots\\
\widehat{f^{\ell_m}}(-s_{m-1}+s_m)&=\sum_{\stackrel{\text{\footnotesize
$(k_{\ell_{m-1}+1},
\ldots ,k_{\ell_m}):$}}
{k_{\ell_{m-1}+1}+\ldots +k_{\ell_m}=s_{m-1}-s_m}}\hat f(k_{\ell_{m-1}+1})
\cdot\ldots\cdot \hat f(k_{\ell_m}),
\end{split}
\end{align*}
we obtain
\begin{align}
\begin{split}
&C_{\ell ,n}(f)=\sum_{k_1+\ldots +k_\ell =0}\hat f(k_1)\cdot\ldots\cdot
\hat f(k_\ell )\cdot\sum^\ell_{m=1}\frac{(-1)^{m-1}}{m}\sum_{\stackrel
{\text{\footnotesize $(\ell_1,\ldots ,\ell_m):$}}{\ell_1+\ldots +\ell_m=
\ell ,\ \ell_i\geq 1}}\\
&\qquad \frac{\ell !}{\ell_1!\ldots\ell_m!}\cdot\#\{u: 0\leq u\leq n-1,
0\leq u+\\
&\qquad \sum^{\ell_1}_1k_i\leq n-1,\ldots ,0\leq u+\sum^{\ell_1+\ldots +
\ell_{m-1}}_1k_i\leq n-1\}.
\end{split}
\end{align}
The last factor in (2.8) is equal to 
\begin{equation}
n-\max \left (0,\sum^{\ell_1}_1k_i,\ldots ,
\sum^{\ell_1+\ldots +\ell_{m-1}}_1k_i \right )-\max \left (0,\sum^{\ell_1}_1(-k_i),
\ldots ,\sum^{\ell_1+\ldots +\ell_{m-1}}_1(-k_i) \right )
\end{equation}
if the expression in (2.9) is nonnegative or zero otherwise.

Lemma 1 is proven.\hfill$\Box$

\section{Proof of the Main Combinatorial Lemma}

First we show that $G(k_1,\ldots ,k_\ell )$ is a linear combination of
terms $\vert k_{i_1}+\ldots +k_{i_s}\vert$.  Then we compute the
coefficient in front of every such term and show it to be equal to zero.

Assume $\ell >2$.
Consider a partition $\p =\{P_1,\ldots ,P_m\}$ of the set $\{1,2,\ldots ,
\ell\}$.  Let us denote $v_1=\sum_{j\in P_1}k_j,\ldots , v_m=\sum_{j\in
P_m}k_j$.  The expression for $G$ can be transformed into
\begin{align}
\begin{split}
&G(k_1,\ldots ,k_\ell )=\sum^\ell_{m=1}\ \ \sum_{\p =\{P_1,\ldots ,P_m\}}
\frac{(-1)^m}{m}\cdot\sum_{\tau\in S_m}\\
&\qquad \max (0,v_{\tau (1)},v_{\tau (1)}+v_{\tau (2)},\ldots ,
v_{\tau (1)}+v_{\tau (2)}+\ldots +v_{\tau (m-1}).
\end{split}
\end{align}

In [R-S] Rudnick and Sarnak, following the ideas of [K] and [Spi] (see 
also [B], [An] ) , used
the following identity for the set of real numbers $v_1,\ldots v_m$ with
zero sum:
\begin{align}
\begin {split}
&\frac{1}{m}\sum_{\tau\in S_m}\max (0,v_{\tau (1)}, 
v_{\tau (1)}+v_{\tau (2)}, \ldots ,v_{\tau (1)}+
v_{\tau (2)}+\ldots v_{\tau (m-1})\\
&\qquad =\frac{1}{4}\sum_{\stackrel{\text{\footnotesize $F\subset 
\{1,\ldots m\},$}}{F,F^C\neq\emptyset}}
(\vert F\vert -1)!(m-\vert F\vert -1)!\cdot\vert\sum_{\ell\in F} v_\ell\vert
\end{split}
\end{align}
The last formula gives us
\begin{align}
\begin{split}
&G(k_1,\ldots ,k_\ell )=\frac{1}{4}\sum^\ell_{m=1}\ \sum_{\p =\{P_1,\ldots
,P_m\}}\ \sum_{\stackrel{\text{\footnotesize $F\subset \{1,\ldots
m\},$}}{F,F^C\neq\emptyset}}
(-1)^{\vert F\vert -1}\cdot \\
&\qquad (\vert F\vert -1)!\cdot\biggl\vert\sum_{i\in\sqcup_{j\in F} P_j}
k_i\biggr\vert\cdot (-1)^{(m-\vert F\vert -1)}\cdot (m-\vert F\vert -1)!.
\end{split}
\end{align}
Le us denote by $A$ the subset $\sqcup_{j\in F}P_j$ of $\{1,2,\ldots ,\ell\}$.
Then $\{P_j\}_{j\in F}$ defines a partition of $A$, and $\{P_j\}_{j\in F^C}$
a partition of $A^C=\{1,2,\ldots ,\ell\}\setminus A$.

We change now the order of summation in (3.3):  first we sum over all nonempty
subsets $A$ of $\{1,2,\ldots ,\ell\}$ and then over all partitions of $A$
and $A^C$:
\begin{align}
\begin{split}
&G(k_1,\ldots ,k_\ell )=\frac{1}{4}\sum_{\stackrel{\text{\footnotesize 
$A\subset\{1,\ldots ,\ell\},$}}{A, A^C\neq\emptyset}}\cdot\biggl (\sum_{\stackrel
{\text{\footnotesize over partitions}}{\u =\{U_1,\ldots ,U_r\}\text{ of }
A}} (-1)^{\vert\u\vert -1}(\vert\u\vert -1)!\biggr )\cdot\\
&\qquad \biggl (\sum_{\stackrel{\text{\footnotesize over partitions}}{\u'
\text{ of }A^C}} (-1)^{\vert\u'\vert -1}\cdot (\vert\u'\vert -1)!\biggr )
\cdot\biggl\vert \sum_{i\in A}k_i\biggr\vert .
\end{split}
\end{align}
Finally we note that
\begin{align*}
\begin{split}
&\sum_{\u=\{U_1,\ldots U_r\}}(-1)^{\vert\u\vert -1}\cdot (\vert\u\vert -1)!
=\sum^{\vert A\vert}_{r=1}\sum_{\stackrel{\text{\footnotesize $
(t_1,\ldots ,t_r):$}}{\sum^r_{i=1}
t_i=\vert A\vert , t_i\geq 1}}(-1)^{r-1}\cdot\\
&\qquad  (r-1)!\frac{\vert A\vert !}{t_1!\cdot\ldots\cdot t_r!}\cdot\frac{1}{r!},
\end{split}
\end{align*}
the expression we already considered in (1.14).  Indeed,
there are exactly $\tfrac{\vert A\vert !}{t_1!\cdot\ldots\cdot
t_r!}\cdot\tfrac{1}{r!}$ different partitions of $A$ such that $\{t_1,
\ldots ,t_r\}=\{\vert U_1\vert ,\ldots ,\vert U_r\vert\}$.  If $\vert A\vert
\geq 2$ this sum is zero.  If $\vert A\vert =1$, then $\vert A^C\vert
=\ell -\vert A\vert\geq 2$ and the second factor in (3.4) equals zero
by the same argument.\hfill$\Box$

Now we turn to a combinatorial lemma first formulated in [Spo].  Let us
denote by $\alpha =(\alpha_1,\ldots ,\alpha_\ell ),\ \beta =(\beta_1,
\ldots ,\beta_\ell )$ vectors with entries $\alpha_j \in \{ 0,1 \}$.  We consider
a lexicographic order on the set of such vectors:  $\alpha <\beta$ iff
$\alpha_j\leq \beta_j,\ j=1,\ldots ,\ell$ and at least for one $j_0\ 
\alpha_{j_0}<\beta_{j_0}$.  Following [Spo] we call such nonzero vectors
branches and a set $T$ of ordered branches $T=\{\alpha^{(1)},\ldots ,
\alpha^{(m)}\},\ \alpha^{(1)}<\alpha^{(2)}<\ldots <\alpha^{(m)},\ 
\vert T\vert =m<\ell$, a tree.  We denote by $T(\ell )$ the set of all
trees formed by a $\ell$-dimensional vectors (branches).  A combinatorial
sum in question is
\begin{equation}
U(k_1,\ldots ,k_\ell )=\sum_{ T\in T(\ell )} (-1)^{\vert T\vert -1}
\cdot\max (0,\alpha\cdot k\vert \alpha\in T )
\end{equation}
Here we used the notation $\alpha\cdot k=\sum^\ell_{j=1}\alpha_j\cdot k_j$.
We call $\max (0,\alpha\cdot k\vert\alpha\in T)$ the maximum of the tree
$T$.  For a warm-up we prove
\medskip

\noindent{\bf Proposition 1}
\begin{align*}
\begin{split}
&U(k_1,\ldots ,k_\ell )+U(-k_1,\ldots ,-k_\ell )=G(k_1,\ldots ,k_\ell,
k_{\ell +1})\\
&\qquad +G(-k_1,\ldots ,-k_\ell, -k_{\ell +1}),
\end{split}
\end{align*}
{\it where }$k_{\ell +1}=-k_1-k_2-\ldots -k_\ell $.
\medskip

\noindent{\bf Remark 7}  Once the proposition is proven
we see of course that $U(k_1,\ldots ,k_\ell )+U(-k_1,\ldots ,-k_\ell )$
is zero for $\ell\geq 2$.

\noindent{\bf Proof of Proposition 1}  
In the above notations
$$
G(k_1,\ldots ,k_\ell, k_{\ell +1})=\sideset{}{'}\sum_{T\in T(\ell +1)}
\frac{(-1)^{\vert T\vert -1}}{\vert T\vert}\cdot\max (0,\alpha\cdot
k'\vert\alpha\in T),
$$
where $k'=(k_1,\ldots ,k_\ell ,k_{\ell +1})$, and the sum $\sum'$ is over
all trees $T\in T(\ell +1)$ such that the largest branch of $T$,
$\alpha^{(\vert T\vert )}$ is less than $D=(1,1,\ldots ,1).$  Similarly,
we can write $U(k_1,\ldots ,k_\ell )=\sum''_{T\in T(\ell +1)}
(-1)^{\vert T\vert -1}\cdot\max (0,\alpha\cdot k'\vert\alpha\in T)$, where
the sum $\sum''$ is over the trees $T\in T(\ell +1)$ such that the
$(\ell +1)^{th}$ coordinate of $\alpha^{(\vert T\vert )}$ is zero.
We define a ``rotation'' on the set of all trees such that $\alpha^{(\vert T\vert 
)}\neq D:  W((\alpha^{(1)},\alpha^{(2)},\ldots ,\alpha^{(\vert T\vert )}))=
(\alpha^{(2)}-\alpha^{(1)},\alpha^{(3)}-\alpha^{(1)},\ldots ,
\alpha^{(\vert T\vert )}-\alpha^{(1)}, D-\alpha^{(1)})$.  Since $
\sum^{\ell +1}_{j=1} k_j=0$, we observe that
\begin{align}
\begin{split}
&\max (0,\alpha\cdot k'\vert\alpha\in T)+\max (0,\alpha\cdot (-k')\vert
\alpha\in T)=\\
&\qquad \max \biggl (0,\alpha\cdot k'\vert\alpha\in W(T)\biggr )
+\max (0,\alpha\cdot (-k')\vert \alpha\in W(T)\biggr ).
\end{split}
\end{align}
The last equality implies
\begin{align*}
\begin{split}
U(k_1,\ldots ,k_\ell )+U(-k_1,\ldots ,-k_\ell )&=\sideset{}{''}
\sum_{T\in T(\ell +1)} (-1)^{\vert T\vert -1}\cdot (\max (0,\alpha\cdot
k'\vert\alpha\in T)\\
&\qquad +\max (0,\alpha\cdot (-k')\vert\alpha\in T))\\
&=\sideset{}{''}\sum_{T\in T(\ell +1)}(-1)^{\vert T\vert -1}\cdot\frac{1}
{\vert T\vert}\cdot\sum^{\vert T\vert -1}_{p=0}\\
&\qquad (\max (0,\alpha\cdot k'\vert\alpha\in W^p(T))\\
&\qquad +\max (0,\alpha\cdot
(-k')\vert\alpha\in W^p(T))\\
& =\sum_{\stackrel{\text{\footnotesize $T\in T (\ell +1)$}}{\alpha^{(\vert T\vert )\neq D}}}
\frac{(-1)^{\vert T\vert -1}}{\vert T\vert}\cdot ( \max (0, \alpha\cdot
k'\vert\alpha\in T)  \\
&+\alpha\cdot (-k')\vert\alpha\in T)   ) \\
&=G(k_1,\ldots ,k_{\ell +1})+G(-k_1,\ldots ,-k_{\ell +1})
\end{split}
\end{align*}
Here we used that for any $T'$ with $\alpha^{(\vert T'\vert )}\neq D$
there exist a unique $T$ with $\alpha^{(\vert T\vert )}_{\ell +1}=0$
and $0\leq p<\vert T\vert$ such that $T'=W^p(T)$.\hfill$\Box$
\medskip

\noindent{\bf Proposition 2}
\begin{equation}
U(k_1,\ldots ,k_\ell)=0\text{ if }\ell\geq 2.
\end{equation}
\phantom{1}\hfill$\Box$

We proceed by induction.

It is easy to check the case $\ell =2$.  Let us assume that the proposition
is true for some $\ell\geq 2$.  Consider $U(k_0,k_1,\ldots ,k_\ell )$.
Since $U$ is a symmetric function we may assume $k_0\leq k_1\leq\ldots
k_\ell$.  The continuity of $U$ implies that it is enough to check (3.7)
for nondegenerate vectors $(k_0,k_1,\ldots ,k_\ell )$.  Therefore we
may assume that the coordinates $k_1,\ldots k_\ell$ are  linearly indepdendent over
the integers.  Fix such $k_1,\ldots k_\ell$ and consider $U$ as a piecewise
linear function of $y=k_0,\ U(y,k)=U(y,k_1,k_2,\ldots ,k_\ell )$.  Our first
claim is that $U(y,k)$ is zero for all negative $y$.  To show this we write
\begin{align*}
\begin{split}
U(y,k)&=\sum_{T\in T(\ell +1)}(-1)^{\vert T\vert -1}\max (0,\alpha\cdot
(y,k)\vert\alpha\in T)\\
&=\sum_{\stackrel{\text{\footnotesize $T\in T(\ell +1):$}}
{\alpha^{(\vert T\vert )}_1=0}}+
\sum_{\stackrel{\text{\footnotesize $T\in T(\ell +1):$}}
{\alpha^{(\vert T\vert )}_1=1,\ 
\alpha^{(\vert T\vert )-1}_1=0}}
+\sum_{\stackrel{\text{\footnotesize $T\in T(\ell +1):$}}
{\alpha^{(\vert T\vert )}_1=1,\ 
\alpha^{(\vert T\vert -1)}_1=1}}
\end{split}
\end{align*}
We denote the three subsums by $U_1, U_2, U_3$.  The first subsum is equal
to
$$\sum_{T\in T(\ell )}(-1)^{\vert T\vert -1}\cdot\max (0,\alpha\cdot
k\vert\alpha\in T),$$
the second --
$$\sum_{T\in T(\ell )}(-1)^{\vert T\vert}\cdot\max (0,\alpha\cdot k\vert
\alpha\in T)+\max (0,y),$$
and by the induction assumptions both are zero.  Now we split the third
subsum in two.  Consider the smallest branch $\alpha\in T$  such that
the first coordinate of $\alpha$ is 1, denote this branch by $\alpha '$
and denote the preceding (may be empty) branch by $\alpha ''$.  We write
$U_3=U_{3,1}+U_{3,2}$, where in $U_{3,1}$ the summation is over $T\in T
(\ell +1)$, such that $\alpha^{(\vert T\vert )}_1=1,\ \alpha^{(\vert T
\vert -1)}_1=1$ and $\alpha '-\alpha ''>(1,0,\ldots ,0)$, and in $U_{3,2}$
the summation is over all other trees from $U_3$.  We establish a one-to-one
correspondence between $U_{3,1}$ and $U_{3,2}$:  for any tree $T_1$ with
$\alpha '-\alpha ''>(1,0,\ldots ,0)$ we construct $T_2=\{\alpha^{(1)},
\ldots ,\alpha '', \alpha ''+(1,0,\ldots ,0), \alpha ',\ldots ,\alpha^{\vert
T\vert}\}$.  Clearly,  $ \vert T_2 \vert = \vert T_1 \vert +1 $, therefore
$$(-1)^{\vert T_1\vert -1}\cdot\max (0,\alpha\cdot (y,k)\vert\alpha\in T_1)=
-(-1)^{\vert T_2\vert -1}\cdot\max (0,\alpha\cdot (y,k)\vert\alpha\in T_2),
$$
and  $U_{3,1}$ and $U_{3,2}$ cancel each other.

Now we assume that $y$ is nonnegative and $0\leq y\leq k_1<k_2<\ldots <k_\ell$.
As we already noted $U(y,k_1,\ldots ,k_\ell )$ is a piecewise linear
continuous function.  We claim that it can change its slope only at $y=0$.
Indeed, $U(y,k_1,\ldots ,k_\ell )$ can change its slope only at the points
of degeneracy of $(y,k_1,\ldots ,k_\ell )$, where $\alpha_0\cdot y
+\alpha\cdot k=\alpha_0'\cdot y+\alpha_0'\cdot k$ and the coordinates of
$(\alpha_0,\alpha),\ (\alpha_0',\alpha')$ take values zero and one.
Because $k$ is a non-degenerate vector we must have
$y+\alpha\cdot k=\alpha'\cdot k$ (or $\alpha\cdot k=y+\alpha'k)$.
Since the tree $T$ contains both branches $(1,\alpha )$ and $(0,\alpha')$
only if $\alpha'\leq\alpha$,  \  the only solution for nonnegative
vector $(y,k)$ must be $y=0,\ \alpha'=\alpha$.  We will finish the
proof of the proposition if we show that $U(y,k)=0$ for sufficiently small positive $y$.
We again write $U=U_1+U_2+U_3$ as before.  Then $U_1=0$ by inductive
assumption and $U_3$ is  zero for sufficiently small $y$ ($U_{3,1}$
and $U_{3,2}$ still cancel each other).  We can write the second subsum
$U_2$ as
\begin{align}
\begin{split}
&\sum_{T\in T(\ell )}(-1)^{\vert T\vert}\cdot\biggl (\max (0,\alpha\cdot k
\vert\alpha\in T)+y\biggr )=\\
&\qquad \sum_{T\in T(\ell )}(-1)^{\vert T\vert}\cdot\biggl 
(\max (0,\alpha\cdot k\vert\alpha\in T) \biggr ) +y\cdot\sum_{T\in T(\ell )}
(-1)^{\vert T\vert}
\end{split}
\end{align}
( the last  sum includes empty tree).  The first term in
(3.8) is zero by inductive assumption and the second is also zero since
$$\sum_{T\in T(\ell )}(-1)^{\vert T\vert}=\sum_{\stackrel{\text{
\footnotesize $\ell_1+\ldots
+\ell_m=\ell +1,$}}{\ell_i\geq 1}}\frac{(-1)^{m-1}}{m}\cdot\frac{(\ell + 1)!}
{\ell_1!\cdot\ldots\cdot\ell_m!}=0.
$$
Proposition 2 is proven.\hfill$\Box$

\section{Orthogonal and symplectic groups.}

We start with the orthogonal case.  The eigenvalues of matrix $M\in 
\text{SO}(2n)$
can be arranged in pairs
$$\exp (i\theta_1),\exp (-i\theta_1),\ldots ,\exp (i\theta_n),\exp
(-i\theta_n), 0\leq\theta_1,\theta_2,\ldots ,\theta_n<\pi .$$
Consider the normalized Haar measure on SO$(2n)$.  The probability
distribution of the eigenvalues is defined by its density (see [We]):
\begin{equation}
P_{\text{SO}2n}(\theta_1,\ldots ,\theta_n)=2\cdot\left (\frac{1}{2\pi}
\right )\cdot\prod_{1\leq i<j<\leq n}(2\cos\theta_i-2\cos\theta_j)^2
\end{equation}
The $k$-point correlation functions are given by (see [So1] )
\begin{equation}
\rho_{n,k}(\theta_1,
\ldots ,\theta_k)=\det \biggl (K^+_{2n-1}(\theta_i,\theta_j)\biggr
)_{1\leq i,j\leq n}
\end{equation}
where
\begin{align}
\begin{split}
&K^+_{2n-1}(x,y)=K_{2n-1}(x,y)+K_{2n-1}(x,-y)=\\
&\qquad \frac{1}{2\pi}\cdot\left (\frac{\sin\left (\frac{(2n-1)(x-y)}{2}\right )}
{\sin\left (\frac{x-y}{2}\right )}
+\frac{\sin\left (\frac{(2n-1)(x+y)}{2}\right )}{\sin\left (\frac{x+y}{2}\right )
}\right ).
\end{split}
\end{align}
In [D-S] and [Jo2] Diaconis-Shahshahani and Johansson studied asymptotic properties
of linear  statistics $\sum^n_{j=1}f(\theta_j)$ where for simplicity we may 
assume that $f$ is real even trigonometric
polynomial, $f(\theta )=\sum^m_{k=1}a_k(\ell^{ik\theta}+\ell^{-ik\theta}),
a_k=\hat f(k), k=1,2,\ldots ,m$.  As before we denote the linear statistics
by $S_n(f)$.  Then  $S_n(f)=\Trace\ (\sum^m_{k=1}
a_k  M^k)$.  It was shown that
\begin{align}
\begin{split}
&E_{2n}\exp\left (t\cdot\sum^n_{j=1}f(\theta_j)\right )=\\
&\qquad \exp \left (t \frac{1}{2} \sum^m_{k=1}\bigl (1+(-1)^k\bigr )
\hat f(k)+\frac{t^2}{2}\sum^m_{k=1}k \hat f(k)^2+\bar 0(1)\right )
\end{split}
\end{align}
which implies the convergence in distribution of $\sum^n_{j=1}f(\theta_j)$ to the
normal law
$$N\left (\frac{1}{2}\cdot\sum^m_{k=1}\left (1+(-1)^k\right )\hat f(k),
\sum^m_{k=1} k\cdot\hat f(k)^2\right ).$$
(Actually (4.4) holds under much weaker conditions --- it is enough
to assume $f\in C^{1+\alpha}([0,\pi ]),\alpha >0  $ ).

{\bf Remark 8 } Similarly to the unitary case (4.4) is equivalent to the 
large $n$ asymptotics result for some determinants, this time Hankel determinants
(see [Jo2], [Jo1]).

Our combinatorial approach allows  to prove  CLT for all $f\in C^1
([0,\pi ])$ as well as to study the local linear statistics $\sum^n_{j=1}
g(L_n\cdot (\theta_j-\theta)), 0<\theta <\pi$.In particular  we establish
\medskip

\noindent{\bf Theorem 2}  {\it Let $g$ be a Schwartz function, $L_n
\rightarrow +\infty, \tfrac{L_n}{n}\rightarrow 0$ and $0<\theta <\pi$.
Then $E_{2n}\sum^n_{j=1}g(L_n\cdot (\theta_j -\theta ))=\tfrac{n}
{L_n\cdot\pi}\cdot\int^\infty_{-\infty}g(x)dx+\bar 0(1)$, and the
centralized random variable $\sum^n_{j=1}g(L_n\cdot (\theta_j-\theta ))
-E_{2n}\sum^n_{j=1}g(L_n\cdot (\theta_j-\theta))$ converges in distribution
to the normal law $N(0,\tfrac{1}{2\pi}\int^\infty_{-\infty}\vert\hat g
(t)\vert^2 \vert t\vert dt)$.
}

Theorem 2 also holds for SO$(2n+1)$ and Sp$(n)$.

Let $M\in \text{SO}(2n+1)$.  Then one of the eigenvalues of $M$ is 1 and
the other $2n$ eigenvalues can be arranged in pairs as before.  The
density of the eigenvalues is equal to
\begin{equation}
P_{\text{SO}(2n+1)}(\theta_1,\ldots ,\theta_n)=\left (\frac{2}{\pi}\right )^2
\cdot\prod_{1\leq i<j\leq n}(2\cos\theta_i-2\cos\theta_j)^2\cdot
\prod^n_{i=1}\sin^2\left (\frac{\theta_i}{2}\right ).
\end{equation}
The formula for the $k$-point correlation function is
\begin{equation}
\rho_{n,k}(\theta_1,\ldots ,\theta_k)=\det\left (K^-_{2n}(\theta_i,\theta_j)
\right )_{i,j=1,\ldots ,k}
\end{equation}
where
\begin{align}
\begin{split}
&K^-_{2n}(x,y)=K_{2n}(x,y)-K_{2n}(x,-y)=\\
&\qquad \frac{1}{2\pi} \left (\frac{\sin\left (n  (x-y)\right )}
{\sin\left (\frac{x-y}{2}\right )}-\frac{\sin (n  (x+y))}{\sin\left (
\frac{x+y}{2}\right )}\right ).
\end{split}
\end{align}
The analogue of (4.4) reads
\begin{align}
\begin{split}
&E_{2n+1}\left (\exp \left (t \sum^n_{j=1}f (\theta_j)\right )\right )=\\
&\qquad \exp\left (t \frac{1}{2}\sum^m_{k=1}\left (-1+(-1)^k\right )
\hat f(k)+\frac{t^2}{2}\sum^m_{k=1}k \hat f(k)^2+\bar 0(1)\right ).
\end{split}
\end{align}
In the symplectic case $M\in\text{Sp}(n)$ the $2n$ eigenvalues again can be
arranged in pairs
$$\exp (i\cdot\theta_i),\exp (-i\cdot\theta_1),\ldots ,\exp (i\cdot\theta_n),
\exp ( -i\cdot\theta_n),0\leq\theta_1,\theta_2,\ldots ,\theta_n<\pi,$$
their density is equal to
\begin{equation}
P_{\text{Sp}(n)}(\theta_1,\ldots ,\theta_n)=\left (\frac{2}{\pi}\right )^n\cdot
\prod_{1\leq i<j\leq n}(2\cos\theta_i-2\cos\theta_j)^2\cdot\prod^n_{i=1}
\sin^2(\theta_i),
\end{equation}
and the formula for $k$-point correlation function is
\begin{equation}
\rho_{n,k}(\theta_1,\ldots ,\theta_k)=\det \left (K^-_{2n+1}(\theta_i,
\theta_j)\right )_{i,j=1,\ldots k.}
\end{equation}
The analogue of (4.4) reads
\begin{align}
\begin{split}
&E_n\left (\exp (t \sum^n_{j=1}f(\theta_j))\right )=\\
&\qquad \exp\left (-t \frac{1}{2}\sum^m_{k=1}\left (1+(-1)^k\right )
\hat f(k)+\frac{t^2}{2}\sum^m_{k=1}k \hat f(k)^2+\bar 0(1).
\right )
\end{split}
\end{align}
We will prove Theorem 2 for SO$(2n)$.  The proofs for SO$(2n+1)$ and Sp$(n)$
are almost identical.
\medskip

\noindent{\bf Proof of Theorem 2}
The arguments from \S 1 imply that it is enough to prove
\medskip

\noindent{\bf Lemma 3}  {\it Let $C_{\ell ,n}(f)$ be the $\ell$-th cumulant
of $\sum^n_{j=1}f(\theta_j),\ell\geq 2$.  Then}
\begin{align}
\begin{split}
&\vert C_{\ell ,n}(f)-\sum_{k_1+\ldots +k_\ell =0}\hat f(k_i)\cdot\ldots
\cdot\hat f(k_\ell) \cdot\frac{1}{2} \biggl (G(k_1,\ldots ,k_\ell )\\
&\qquad +G(-k_1,\ldots ,-k_\ell)\biggr )\biggl\vert\leq \const_\ell  
\sum_{\stackrel{\text{\footnotesize $k_1+\ldots +k_\ell =0$}}{\vert k_1\vert +\ldots +\vert k_\ell\vert >n}}\\
&\qquad \vert k_1 \vert \vert\hat f(k_1)\vert\cdot\ldots\cdot \vert \hat f(k_\ell ) \vert +
\const'_\ell \sum_{\vert k_1\vert +\ldots +\vert k_\ell\vert >n}
\vert\hat f(k_1)\vert\cdot\ldots\cdot\vert\hat f(k_\ell )\vert
\end{split}
\end{align}
\hfill$\Box$

We start with the formula (2.6) which holds for general determinantal random
point fields:
\begin{align*}
\begin{split}
&C_{\ell ,n}(f)=\sum^\ell_{m=1}\sum_{\stackrel{\text{\footnotesize
$\ell_1+\ldots +\ell_m=\ell ,$}}{\ell_i\geq 1}}(-1)^m\frac{\ell !}
{\ell_1!\cdot\ldots\cdot\ell_m!}\cdot\frac{1}{m}
\int_{[0,\pi ]^m}f^{\ell_1}(x_1)\cdot\ldots\cdot f^{\ell_m}
(x_m)\cdot\\
&\qquad \prod^m_{j=1} \left ( K_{2n-1}(x_j,x_{j+1})+K_{2n-1}(x_j,-x_{j+1}) \right )
dx_1\ldots dx_m
\end{split}
\end{align*}
(we always assume $x_{m+1}=x_1$).
\begin{align}
\begin{split}
&=\sum^\ell_{m=1}\sum_{\stackrel{\text{$\ell_1+\ldots +\ell_m=\ell
$}}{\ell_i\geq 1}}(-1)^m\cdot\frac{1}{m}\cdot\frac{\ell !}{\ell_1!
\cdot\ldots\cdot\ell_m!}\sum_{\epsilon_1=\pm 1}\sum_{\epsilon_2=\pm 1}\ldots
\sum_{\epsilon_m=\pm 1}\\
&\quad \int_{[0,\pi ]^m}f^{\ell_1}(x_1)\cdot\ldots\cdot f^{\ell_m}(x_m)
\cdot\prod^m_{j=1}K_{2n-1}(x_j,\epsilon_j\cdot x_{j+1})dx_1\cdot\ldots
\cdot dx_m
\end{split}
\end{align}
Each term in the last sum with $\prod^m_{i=1}\epsilon_i=1$ is equal to
\begin{align*}
\begin{split}
&\int_{\prod^m_{i=1}\epsilon_{i-1}\cdot [0,\pi ]}
f^{\ell_1}(x_1)\cdot
\ldots\cdot f^{\ell_m}(x_m)\cdot\prod^m_{j=1}K_{2n-1}(x_j,x_{j+1})
\prod^m_{i=1}d(\epsilon_{i-1}\cdot x_i)=\\
&\qquad \frac{1}{2^m}\cdot\int_{[0,2\pi ]^m}f^{\ell_1}(x_1)\cdot
\ldots\cdot f^{\ell_m}(x_m)\cdot\prod^m_{j=1}K_{2n-1}(x_j,x_{j+1})dx_1
\cdot\ldots\cdot dx_m
\end{split}
\end{align*}
(we use the fact that $f(x)$ is even). Combining these  terms
together we obtain the same expression as for $\tfrac{1}{2}\cdot C_{\ell ,
2n-1}(\sum^{2n-1}_{j=1}f(\theta_j))$ in the case of $U(2n-1)$, 
which gives vanishing contribution if $ \ell > 2 $.  Finally
we claim that the contribution from the terms with $\prod^m_{i=1}
\epsilon_i=-1$ can be bounded from above by
$$\const'_\ell\cdot\sum_{\vert k_1\vert +\ldots +\vert k_\ell\vert
>n}\vert\hat f(k_1)\vert\cdot\ldots\cdot\vert\hat f(k_\ell )\vert .
$$
Indeed, the integral
$$\int_{[0,\pi ]^m}f^{\ell_1}(x_1)\cdot\ldots\cdot f^{\ell_m}(x_m)\cdot
\prod^m_{j=1}\left (\frac{1}{2\pi}\sum^n_{s_j=-n} e^{is_j(x_j-\epsilon_j
\cdot x_{j+1})}\right )dx_1\ldots dx_m$$
can be rewritten as
\begin{align*}
\begin{split}
\frac{1}{2^m}\cdot\sum^n_{s_1=-n}\ldots\sum^n_{s_m=-n}\widehat{f^{\ell_1}}
(s_1-\epsilon_m\cdot s_m)\cdot\widehat {f^{\ell_2}}(s_2-\epsilon_1\cdot
s_1)\cdot\ldots\cdot \widehat {f^{\ell_m}}(s_m-\epsilon_{m-1}\cdot
s_{m-1})
\end{split}
\end{align*}
Consider the euclidian basis $\{e_j\}^m_{j=1}$ in $\bbR^m$ and define $f_j=e_j
-\epsilon_{j-1} e_{j-1},\ \epsilon_0=\epsilon_m$.  The vectors
$\{f_j\}^m_{j=1}$ form a basis in $\bbR^m$ iff $\prod^m_{j=1}\epsilon_j=-1$.
Then for any $m$-tuple $(t_1,\ldots ,t_m)$ there exists the only $m$-tuple
$(s_1,\ldots ,s_m)$ such that $t_j=s_j-\epsilon_{j-1}\cdot s_{j-1},\ 
j=1,\ldots ,m$.  
We write $\widehat {f^{\ell_j}}(t_j)=\sum\hat f(k_{\ell_1+\ldots +
\ell_{j-1}+1})\cdot\ldots\cdot\hat f(k_{\ell_1+\ldots +\ell_j})$, where
the sum is over $k_i$ such that $\sum^{\ell_1+\ldots +\ell_j}_{\ell_1+
\ldots +\ell_{j-1}+1}k_i=t_j$.  When we plug this into (4.13) we obtain
a linear combination of
\begin{equation}
\hat f(k_1)\cdot\ldots\cdot \hat f(k_m)
\end{equation}
It is easy to see that for $\vert k_1\vert +\ldots +\vert k_m\vert\leq n$
the coefficient with the term (4.14) is equal to
$$\frac{1}{2^m}\cdot\sum^\ell_{m=1}\ \sum_{\stackrel{\text{\footnotesize
$\ell_1+\ldots +\ell_m=\ell ,$}}{\ell_i\geq 1}}(-1)^m\cdot\frac{1}{m}
\cdot\frac{\ell !}{\ell_1!\cdot\ldots\cdot\ell_m!}=0
$$
For $\vert k_1\vert +\ldots +\vert k_m\vert >n$ the coefficient is bounded 
from above by some constant.  This finished the proof of Lemma 3.
\hfill$\Box$

Similar to \S 1 we obtain the proof of Theorem 2 by applying the lemma to
$\sum^n_{j=1}g(L_n\cdot (\theta_j-\theta ))$.

\def\cmp{{\it Commun. Math. Phys., }}
\def\jsp{{\it J. Stat. Phys., }}
\def\dmj{{\it Duke Math. J., }}
\def\am{{\it Ann. of Math., }}


\begin{thebibliography}{ZZZZZ}
\bibitem[A]{}
C. Andr\'eief, Note sur une relation les int\'egrales d\'efinies des
produits des fonctions, {\it M\'em. de la Soc. Sci. Bordeaux},
{\bf 2}, 1--14, 1883.


\bibitem[An]{}
E.S.Andersen, On sums of symmetrically dependent random variables,
{\it Skand. Aktuarietidskr.},
{\bf 36}, 123-138, 1953.


\bibitem[B]{}
G.Baxter,  Combinatorial methods in fluctuation theory,
{\it Z.Wahrscheinlichkeitstheorie},
{\bf 1}, 263-270, 1963.


\bibitem[Ba]{}
E. Basor, Distribution functions for random variables for ensembles of
positive Hermitian matrices, {\it Comm. Math. Phys. }, {\bf 188}, 327--350, 1997.


\bibitem[B-F]{}
T. H. Baker, and P.J. Forrester, Finite $N$ fluctuation formulas
for random matrices, \jsp {\bf 88}, 1371--1385, 1997.


\bibitem[Ba-W]{}
E. Basor and H.Widom,Toeplitz and Wiener-Hopf determinants with piecewise 
contionuous symbols,
{\it
J. of Funct. Anal.}, {\bf 50},387-413, 1983


\bibitem[BM-K]{}
A. Boutet de Monvel and A. Khorunzky, Asymptotic distribution of
smoothed eigenvalue density, I, II, {\it Rand. Oper.
Stoch. Eqn}, {\bf vol.7, No. 1}, 1-22 , 1999  and {\bf vol.7, No.2},
149-168, 1999.


\bibitem[Bo]{}
A.B\"ottcher, The Onsager formula, the Fisher-Hartwig conjecture, 
and their influence on research into Toeplitz operators,
 \jsp
 {\bf 78}, 575-588 , 1995.


\bibitem[Bo-S]{}
A. B\"ottcher and B.Silbermann,
{\it Introduction to Large Truncated Toeplitz matrices },
Springer, 1999.


\bibitem[C-L]{}
O. Costin and J. Lebowitz, Gaussian fluctuations in random matrices, 
{\it Phys. Rev. Lett.}, {\bf 75}, 69--72, 1995.

\bibitem[De]{}
A.Devinatz, The strong Szeg\"o limit theorem,
{\it Illinois J. Math.}, {\bf 11}, 160-175, 1967.


\bibitem[D-S]{}
P. Diaconis and M. Shahshahani, On the eigenvalues of random matrices,
Studies in Appl. Probab., Essays in honor of Lajos Takacs, {\it J.
Appl. Probab.}, {\bf Special Volume 31A}, 49--62, 1994.

\bibitem[D]{}
P. Diaconis , Patterns in eigenvalues,
{\it  Bull. Amer. Math. Soc. },
to appear.


\bibitem[Dy]{}
F.J. Dyson, Correlations between eigenvalues of a random matrix,
{\it  Comm. Math. Phys.}, {\bf 19}, 235-250, 1970.

\bibitem[F-H]{}
H.E. Fisher and R.E. Hartwig, Toeplitz determinants, some applications,
theorems and conjectures, {\it Adv. Chem. Phys.}, {\bf 15}, 333--353, 
1968.

\bibitem[H]{}
I.I. Hirschman, Jr.,
On a theorm of Szeg\"o, Kac and Baxter, {\it J. d' Analyse Math.}, {\bf 14}, 225-234,
1965.

\bibitem[G-I]{}
B.L. Golinskii and I.A. Ibragimov, On Szeg\"o's limit theorem,
{\it Math. USSR-Izv}, {\bf vol 5, No. 2}, 421--446, 1971.

\bibitem[I-D]{}
C. Itzykson and J-M. Drouffe, {\it Statistical Field Theory}, vol. 1,
Cambridge University Press, 1989.

\bibitem[Jo1]{}
K. Johansson, On Szeg\"o's asymptotic formula for Toeplitz determinants
and generalizations, {\it Bull. Sci. Math.}, {\bf 112}, 257--304,
1988.

\bibitem[Jo2]{}
K. Johansson, On random matrices from classical compact groups, \am {\bf
145}, 519--545, 1997.

\bibitem[Jo3]{}
K. Johansson, On fluctuation of eigenvalues of random Hermitian matrices,
\dmj {\bf 91}, 151--204, 1998.

\bibitem[K]{}
M. Kac, Toeplitz matrices, translation kernels and a related problem
in probability theory, \dmj {\bf 21}, 501--509, 1954.

\bibitem[KKP]{}
A. Khorunzhy, B. Khoruzhenko, and L. Pastur, Asymptotic properties of
large random matrices with independent entries, {\it J. Math. Phys.}, 
{\bf 37}, 5033--5059, 1996.

\bibitem[MC-W]{}
B.M. McCoy and T.T. Wu, {\it The Two-dimensional Ising Model}, Harvard
University Press, Cambridge, Massachusetts, 1973.

\bibitem[Me]{}
M.L. Mehta, {\it Random Matrices}, 2nd edition, Academic Press, Boston,
1991.

\bibitem[R-S]{}
Z. Rudnick and P. Sarnak, Zeroes of principal $L$-functions and random
matrix theory, A Celebration of John F. Nash, Jr., \dmj {\bf 61},
269--322, 1996.

\bibitem[SSo1]{}
Ya. Sinai and A. Soshnikov, A refinement of Wigner's semicircle law in a
neighborhood of the spectrum edge for random symmetric matrices,
{\it Funct. Anal. Appl.}, {\bf vol 32, no. 2}, 114--131, 1998.

\bibitem[SSo2]{}
Ya. Sinai and A. Soshnikov, Central limit theorem for traces of
large random matrices with independent entries, {\it Bol. Soc.
Brasil. Mat.}, {\bf vol. 29, no. 1}, 1--24, 1998.

\bibitem[So1]{}
A. Soshnikov, Level spacings distribution for large random matrices: 
gaussian fluctuations, \am {\bf 148}, 573--617, 1998.

\bibitem[So2]{}
A. Soshnikov, Gaussian fluctuations in Airy, Bessel, sine and other
determinantal random point fields,
Preprint, 1999, available via http://xxx.lanl.gov/abs/math/9907012.

\bibitem[Spi]{}
F. Spitzer, A combinatorial lemma and its applications to probability
theory, {\it Trans. Amer. Math. Soc.}, {\bf 82}, 323--339, 1956.

\bibitem[Spo]{}
H. Spohn, Interacting Brownian particles:  A study of Dyson's model,
in {\it Hydrodynamic Behavior and Interacting Particle Systems},
G. Papanicolau, ed., Springer-Verlag, New York, 1987.

\bibitem[Sz]{}
G. Szeg\"o, On certain Hermitian forms associated with the Fourier
series of a positive function, {\it Comm. Seminaire Math. de l'Univ.
de Lund, tome suppl\'ementaire, d\'edi\'e  \'a Marcel Riesz}, 228--237,
1952.

\bibitem[T-W]{}
C.A. Tracy and H. Widom, Correlation functions, cluster functions,
and spacing distributions for random matrices, {\it J. Stat. Phys.},
{\bf vol. 92, no. 5/6}, 809--835, 1998.

\bibitem[We]{}
H. Weyl, {\it The Classical Groups:  Their Invariants and Representations},
Princeton Univ. Press, Princeton, 1939.

\bibitem[Wid1]{}
H. Widom, Toeplitz determinants with single generating function,
{\it
Am.  J. Math.}, {\bf 95}, 333-383, 1973.

\bibitem[Wid2]{}
H. Widom, Asymptotoc behaviour of block Toeplitz matrices and determinants,
I and II,
{\it
Adv. Math.}, {\bf 13}, 284-322, 1973 and {\bf 21}, 1-29, 1976.

\bibitem[Wie]{}
K. Wieand, Eigenvalue distributions of random matrices in the
permutation group and compact Lie groups, Ph.D. thesis, Dept. Math,
Harvard, 1998.


\end{thebibliography}
\end{document}